\documentclass[english,11pt]{article}
\usepackage[latin1]{inputenc}
\usepackage{amsmath,amsthm,amssymb,babel}

\newtheorem{teo}{Theorem}
\newtheorem{defin}{Definition}

\newtheorem{lemma}{Lemma}

\def\proof{\noindent{\it Proof.}\ }
\def\endproof{\hfill $\Box$\par\vskip3mm}

\def\eq#1{(\ref{#1})}

\def\neweq#1{\begin{equation}\label{#1}}
\def\endeq{\end{equation}}
\def\weak{\rightharpoonup}

\def\de{\partial}
\def\phi{\varphi}
\def\RR{{\mathbb R} }

\def\ri{\rightarrow}
\def\intom{\int_\Omega}
\date{}

\title{\sc On a class of variational-hemivariational inequalities
involving set valued mappings }
\author{\sc Nicu\c sor Costea$\;^{a,}$\thanks{Correspondence address:
Nicu\c sor Costea, Department of Mathematics, University of
Craiova, 200585 Craiova, Romania. E-mail: {\tt
nicusorcostea@yahoo.com}} \  and Cezar Lupu$\;^{a,b}$\\
\small $^a$Department of Mathematics, University of Craiova,
200585 Craiova,
Romania\\
\small $^b$Department of Mathematics, University of Bucharest, 70109 Bucharest, Romania\\
\small E-mail addresses: {\tt nicusorcostea@yahoo.com} \qquad {\tt
lupucezar@gmail.com}}

\begin{document}
\baselineskip16pt \maketitle

\noindent{\small{\sc Abstract}. Using the KKM technique,  we
establish some existence results for variational-hemivariational
inequalities involving monotone set valued mappings on bounded,
closed and convex subsets in reflexive Banach spaces. We also
derive several sufficient
conditions for the existence of solutions in the case of unbounded subsets.
}

\section{Introduction}

Because of its wide applicability in various domains such as
mechanics, engineering sciences, economics, optimal control,
inequality theory has become an important branch of mathematics.
Inequality problems can be divided into two main classes: that of
variational inequalities and that of hemivariational inequalities.
The study of variational inequalities began with the work of G.
Fichera \cite{Fich}, J. L. Lions and G. Stampacchia
\cite{Lio-Sta}. A basic result in the theory of variational
inequalities is due to Hartman and Stampacchia \cite{Har-Sta} and
states that if $X$ is a finite dimensional Banach space, $K\subset
X$ is compact and convex, and $A$ is a continuous operator, then
the variational inequality problem of finding $u\in K$ such that
\neweq{VI}
\langle Au,v-u\rangle\geq 0,\quad\forall v\in K
\endeq
has a solution. When $K$ is not compact, or $X$ is infinite
dimensional, certain monotonicity properties are required to prove
the existence of solution.

\medskip

In 1999, Panagiotopoulos, Fundo and R\u adulescu
\cite{Pan-Fun-Rad} extended the classical results from
\cite{Har-Sta}, proving several versions of theorems of
Hartman-Stampacchia's type for the case of hemivariational
inequalities on compact or on closed and convex subsets, in
infinite and finite dimensional Banach spaces.
\smallskip

By replacing the subdifferential of a convex function by the
generalized gradient (in the sense of F. H. Clarke) of a locally
Lipschitz functional, hemivariational inequalities arise whenever
the energetic functional associated to a concrete problem is
nonconvex. This new type of inequalities appears as a
generalization of the variational inequalities, but
hemivariational inequalities are much more general, in the sense
that they are not equivalent to minimum problems but, give rise to
substationarity problems. The theory of hemivariational
inequalities can be viewed as a new field of Nonsmooth Mechanics
since the main ingredient used in the study of these inequalities
is the notion of Clarke subdifferential of a locally Lipschitz
functional. The mathematical theory hemivariational inequalities,
as well as their applications in Mechanics, Engineering or
Economics were introduced and developed by P. D. Panagiotopoulos
\cite{Pan1}-\cite{Pan5} in the case of nonconvex energy functions.
For a treatment of this theory and further comments we recommend
the monographs by Z. Naniewicz and P. D. Panagiotopoulos
\cite{Nan-Pan}, D. Motreanu and P. D. Panagiotopoulos
\cite{Mot-Pan} and by D. Motreanu and V. R\u adulescu
\cite{Mot-Rad2}. For more information and connections regarding
hemivariational inequalities see the recent papers of Costea-R\u
adulescu \cite{Cos-Rad1,Cos-Rad2}, Mig\' orski-Ochal-Sofonea
\cite{MOG1,MOG2}, Z. Liu-G. Liu \cite{Liu} and the references
therein.

\section{The abstract framework}
Throughout this paper $X$ will denote a real reflexive Banach
space with its dual space $X^\ast$ and $T:X\ri L^p(\Omega;\RR^k)$
will  be a linear and compact operator where $1<p<\infty$ and
$\Omega$ is a bounded and open subset of $\RR^N$. We shall denote
$Tu=\hat u$ and by $p'$ the conjugated exponent of $p$. Let
$j=j(x,y):\Omega\times\RR^k\ri \RR$ be a Carath\' eodory function,
locally Lipschitz with respect to the second variable which
satisfies the following condition:

-there exist $C>0$ such that
\neweq{eq1}
|z|\leq C(1+|y|^{p-1})
\endeq
a.e. $x\in \Omega$, for all $y\in \RR^k$ and all $z\in\de j(x,y)$.

We shall use the notation $j^{\,0}(x,y;h)$ for the Clarke's
generalized directional derivative (see e.g. \cite{Cla} or
\cite{Mot-Rad2}) of the locally Lipschitz mapping $j(x,\cdot)$ at
the point $y\in\RR^k$ with respect to the direction $h\in\RR^k$,
where $x\in\Omega$, i.e.,
$$
j^{\,0}(x,y;h)=\limsup\limits_{\stackrel{w\ri y}{\lambda\downarrow
0}}\frac{j(x,w+\lambda h)-j(x,w)}{\lambda}.
$$
Accordingly, Clarke's generalized gradient $\de j(x,y)$ of the
locally Lipschitz map $j(x,\cdot)$ is defined by
$$
\de j(x,y)=\{z\in \RR^k|\  z\cdot h  \leq j^{\,0}(x,y;h), \mbox{
for all } h\in\RR^k\},
$$
where the symbol ``$\cdot$'' means the inner product on $\RR^k$.

The euclidian norm in $\RR^k$, $k\geq 1$, resp. the duality
pairing between a Banach space and its dual will be denoted by
$|\cdot|$, respectively $\langle\cdot,\cdot\rangle$. We also
denote by $\|\cdot\|_p$ the norm in the space $L^p(\Omega;\RR^k)$
defined by
$$
\|\hat u\|_p=\left(\intom |\hat u(x) |^p \; dx\right)^{1/p}.
$$
Let $K$ be a nonempty closed, convex subset of $X$ and
$\phi:X\ri(-\infty,+\infty]$ a convex lower semicontinuous
functional such that
\neweq{eq3}
D(\phi)\cap K\neq\emptyset.
\endeq

Assuming $A$ is a set valued mapping from $K$ into $X^\ast$, i.e.
precisely into the power set $\mathcal{P}(X^\ast)$ or $2^{X^\ast}$
of $X^\ast$, with the domain of $A$, i.e. the set $D(A)=\{u\in K|\
A(u)\neq\emptyset\}$, the whole set $K$ our aim is to study the
following variational-hemivariational inequality problem:

Find $u\in K$ and $u^\ast\in A(u)$ such that
\neweq{eq4}
\langle u^\ast,v-u\rangle+\phi(v)-\phi(u)+\intom j^0(x,\hat
u(x);\hat v(x)-\hat u(x))\; dx\geq 0,\quad \forall v\in K.
\endeq

\noindent As it will be seen, this problem closely links to the
following problem:

Find $u\in K$ such that
\neweq{eq5}
\sup\limits_{v^\ast\in A(v)} \langle v^\ast,u-v\rangle\leq
\phi(v)-\phi(u)+\intom j^0(x,\hat u(x);\hat v(x)-\hat u(x))\;
dx,\quad\forall v\in K.
\endeq
\begin{itemize}
\item We remark that if $T\equiv 0$ and $\phi$ is the indicator
function of the set $K$, that is
$$
\phi\equiv I_K=\left\{
\begin{array}{ll}
0,& \mbox { if } u\in K\\
\infty,& \mbox{ otherwise},
\end{array}
\right.
$$
then problem \eq{eq4} is equivalent to finding $u\in K$ and
$u^\ast\in A(u)$ such that
\neweq{eqSC1}
\langle u^\ast,v-u\rangle\geq0,\quad\forall v\in K,
\endeq
which is called the {\it generalized variational inequality
problem}, while problem \eq{eq5} becomes equivalent to finding
$u\in K$ such that
\neweq{eqSC2}
\sup\limits_{v^\ast\in A(v)}\langle v^\ast,u-v\rangle\leq
0,\quad\forall v\in K
\endeq
which is called the {\it dual variational inequality problem} and
was originally discussed in \cite{Min} ( see also \cite{Bro,
Kon}).

\item We note that if $A:K\ri X^\ast$ is a single valued operator
and $\phi\equiv I_K$ problem \eq{eq4} is equivalent to finding
$u\in K$ such that
\neweq{eqSC3}
\langle Au,v-u\rangle+\intom j^0(x,\hat u(x);\hat v(x)-\hat
u(x))\geq 0,\quad\forall v\in K.
\endeq
The inequality of the type \eq{eqSC3} is known as the {\it
hemivariational inequality of Hartman-Stampachhia type}. It turned
out that various problems from Nonsmooth Mechanics can be studied
using this inequality, see \cite{Pan-Fun-Rad}.

\item If $A$ is a single valued and $T\equiv 0$ then problem
\eq{eq4} becomes equivalent to finding $u\in K$ such that
\neweq{eqSC4}
\langle Au,v-u\rangle+\phi(v)-\phi(u)\geq0,\quad\forall v\in K
\endeq
which is known as the {\it mixed variational inequality}. For
applications, numerical methods and formulations, see
\cite{Bai-Cap,Gia-Mau,Glo-Lio-Tre,Noor}.

\item Finally, if $A$ is single valued, $T\equiv 0$ and
$\phi\equiv I_K$ then problem \eq{eq4}  becomes the {\it standard
variational inequality}: find $u\in K$ such that
\neweq{eqSC5}
\langle Au,v-u\rangle\geq0,\quad\forall v\in K.
\endeq
\end{itemize}

It is clear that problems \eq{eqSC1}, \eq{eqSC3}-\eq{eqSC5} are
special cases of our variational-hemivariational inequality
problem \eq{eq4}. In brief, for a suitable choice of $A$, $\phi$
and $T$ one can obtain a wide class of inequality problems. This
shows that \eq{eq4} is quite general and unifying in the same
time.

\begin{defin}\label{defin1}
A set valued mapping $A:K\ri 2^{X^\ast}$ is said to be
\begin{description}
\item (i) monotone on $K$ if, for each pair of point $u,v\in K$
and all $u^\ast\in A(u)$ and $v^\ast\in A(v)$,
$$
\langle v^\ast-u^\ast, v-u\rangle\geq 0;
$$
\item (ii) lower semicontinuous  at $u\in K$ (see e.g.
\cite{Nan-Pan}, p. 12), if and only if for any $u^\ast \in A(u)$
and for any sequence $\{u_n\}\subset K$ with $u_n\ri u$, a
sequence $u_n^\ast\in A(u_n)$ can be determined which converges to
$u^\ast$;

\item(iii) lower hemicontinuous on $K$ if the restriction of $A$
to every line segment of $K$ is lower semicontinuous with respect
to the weak-star topology in $X^\ast$.
\end{description}
\end{defin}

\section{Existence results}

In this section, we discuss the existence of solutions of the
variational-hemivariational inequality problem \eq{eq4}. The
following two results will be needed in the sequel.

\begin{lemma}\label{lema1}(see \cite{Fan} or \cite{KKM})
Let $K$ be a nonempty subset of a Hausdorff topological vector
space $E$ and let $G:K\ri 2^E$ be a set valued mapping satisfying
the following properties:
\begin{description}
\item (i) $G$ is a KKM mapping;

\item (ii) $G(x)$ is closed in $E$ for every $x\in K$;

\item (iii) $G(x_0)$ is compact in $E$ for some $x_0\in K$.
\end{description}
Then $\bigcap_{x\in K} G(x)\neq \emptyset$.
\end{lemma}
We recall that a set valued mapping $G:K\ri 2^E$ is said to be a
{\it KKM mapping} if for any $\{x_1,\ldots,x_n\}\subset K$,
$\rm{co}\{x_1,\ldots,x_n\}\subset \bigcup_{i=1}^n G(x_i)$, where
$\rm{co}\{x_1,\ldots,x_n\}$ denotes the convex hull of
$x_1,\ldots,x_n$.

\begin{lemma}\label{lema2}
\cite{Pan-Fun-Rad} If $j$ satisfies \eq{eq1}, $X_1,X_2$ are
nonempty subsets of $X$ and $T:X\ri L^p(\Omega;\RR^k)$ is  a
linear and compact operator, then the mapping $X_1\times
X_2\ri\RR$ defined by
$$
(u,v)\longmapsto\intom j^0(x,\hat u(x);\hat v(x)-\hat u(x))\; dx
$$
is weakly upper semicontinuous.
\end{lemma}
\smallskip

We denote by $S$ and $S^\star$ the solutions sets of problem
\eq{eq4} and problem \eq{eq5}, respectively. Now we give the
relationship between problem \eq{eq4} and \eq{eq5}.
\begin{teo}\label{Th1}
Let $K$ be a nonempty closed and convex subset of the real
reflexive Banach space $X$. If $A:K\ri 2^{X^\ast}$ is monotone,
then $S\subseteq S^\star$. In addition, if $A$ is lower
hemicontinuous, then $S^\star=S$.
\end{teo}

\proof Let $u\in S$ and $v\in K$ be arbitrary fixed. Then it
exists $u^\ast\in A(u)$ such that
\neweq{eq6}
\langle u^\ast,v-u\rangle+\phi(v)-\phi(u)+\intom j^0(x,\hat
u(x);\hat v(x)-\hat u(x))\; dx \geq0.
\endeq
Since $A$ is monotone we have
\neweq{eq7}
\langle v^\ast-u^\ast,v-u\rangle\geq 0,\quad \forall v^\ast\in
A(v)
\endeq
Hence by \eq{eq6} and \eq{eq7} we have
\neweq{eq8}
\langle v^\ast,v-u\rangle+\phi(v)-\phi(u)+\intom j^0(x,\hat
u(x);\hat v(x)-\hat u(x))\; dx \geq0,\quad\forall v^\ast\in A(v).
\endeq
This is equivalent to
\neweq{eq9}
\sup\limits_{v^\ast\in A(v)}\langle
v^\ast,u-v\rangle\leq\phi(v)-\phi(u)+\intom j^0(x,\hat u(x);\hat
v(x)-\hat u(x))\; dx.
\endeq
Since $v$ has been arbitrary chosen, it follows that \eq{eq9}
holds for all $v\in K$ which implies that $u\in S^\star$.
\medskip

In addition if $A$ is lower hemicontinuous, we will show that
$S=S^\star$. Suppose $u\in S^\star$ and let $v\in K$ be arbitrary
fixed. We define the sequence $\{u_n\}_{n\geq1}$ by
$u_n=u+\frac{1}{n}(v-u)$. Clearly $\{u_n\}\subset K$ by the
convexity of $K$. For any $u^\ast\in A(u)$, using the lower
hemicontinuity of $A$, a sequence $u_n^\ast\in A(u_n)$ can be
determined which converges weakly-star to $u^\ast$. Taking into
account that $u\in S^\star$, for each $n\geq 1$ we have
\neweq{eq10}
\langle u_n^\ast,u-u_n\rangle\leq \phi(u_n)-\phi(u)+\intom
j^0(x,\hat u(x);\hat u_n(x)-\hat u(x))\; dx
\endeq
But $\phi$ is convex and $j^0(x,\hat u;\lambda\hat v)=\lambda
j^0(x,\hat u;\hat v)$ for all $\lambda>0$. Therefore \eq{eq10} may
be written, equivalently
\begin{eqnarray*}
0&\leq& \left\langle
u_n^\ast,\frac{1}{n}(v-u)\right\rangle+\phi\left(\frac{1}{n}v+
\frac{n-1}{n}u\right)-\phi(u)+\intom
j^0\left(x,\hat u(x);\frac{1}{n}(\hat v(x)-\hat u(x))\right)\; dx\\
&\leq&\frac{1}{n}\left[\langle
u_n^\ast,v-u\rangle+\phi(v)-\phi(u)+\intom j^0(x,\hat u(x);\hat
v(x)-\hat u(x))\; dx\right].
\end{eqnarray*}
Multiplying the last relation by $n$ and passing to the limits as
$n\ri\infty$ we obtain the $u\in S$.
\endproof

\begin{teo}\label{Th2}
Let $K$ be a nonempty, bounded, closed and convex subset of the
real reflexive Banach space $X$ and $A:K\ri 2^{X^\ast}$ a set
valued mapping which is monotone and lower hemicontinuous on $K$.
If $T:X\ri L^p(\Omega;\RR^k)$ is linear and compact and $j$
satisfies the condition \eq{eq1} then the
variational-hemivariational inequality problem \eq{eq4} has at
least one solution.
\end{teo}

\proof For any $v\in K\cap D(\phi)$ define two set valued mappings
$F,G:K\cap D(\phi)\ri 2^X$ as follows:
$$
F(v)=\left\{u\in K\cap D(\phi)\left|
\begin{array}{l}
 \mbox{ there exists }
u^\ast\in A(u)\mbox{ such that }\\
\langle u^\ast,v-u\rangle+\phi(v)-\phi(u)+\intom j^0(x,\hat
u(x);\hat v(x)-\hat u(x))\; dx\geq 0
\end{array}
\right\}\right.
$$
and
$$
G(v)=\left\{u\in K\cap D(\phi)\left|\sup\limits_{v^\ast\in
A(v)}\langle v^\ast,u-v\rangle\leq \phi(v)-\phi(u)+\intom
j^0(x,\hat u(x);\hat v(x)-\hat u(x))\; dx\right.\right\}.
$$
\begin{description}
\item {\sc Step 1.} $F$ is a KKM mapping.

If $F$ is not a KKM mapping, then there exists
$\{v_1,\ldots,v_n\}\subset K\cap D(\phi)$ such that
$$
\rm{co}\{v_1,\ldots,v_n\}\not\subset\bigcup\limits_{i=1}^n F(v_i)
$$
i.e., there exists a $v_0\in \rm{co}\{v_1,\ldots,v_n\}$,
$v_0=\sum_{i=1}^n \lambda_i v_i$, where $\lambda_i\in [0,1]$,
$i=1,\ldots,n$, $\sum_{i=1}^n \lambda_i=1$, but $v_0\not\in
\bigcup_{i=1}^n F(v_i)$. By the definition of $F$, we have
$$
\langle v_0^\ast, v_i-v_0\rangle+\phi(v_i)-\phi(v_0)+\intom
j^0(x,\hat v_0(x);\hat v_i(x)-\hat v_0(x))\; dx<0,\quad\forall
v_0^\ast\in A(v_0)
$$
for $i=1,\ldots,n$. It follows from the convexity of the mapping
$\hat v\longmapsto j^0(x,\hat u;\hat v)$ and the convexity of
$\phi$ that for each  $v_0^\ast\in A(v_0)$ we have
\begin{eqnarray*}
0&=&\langle v_0^\ast, v_0-v_0\rangle+\phi(v_0)-\phi(v_0)+\intom
j^0(x,\hat v_0(x);\hat v_0(x)-\hat v_0(x))\; dx\\
&=&\left\langle v_0^\ast,\sum_{i=1}^n \lambda_i
v_i-v_0\right\rangle+\phi\left(\sum_{i=1}^n \lambda_i
v_i\right)-\phi(v_0)+\intom j^0\left(x,\hat
v_0(x);\sum_{i=1}^n\lambda_i\hat v_i(x)-\hat v_0(x)\right)\;dx\\
&\leq&\sum_{i=1}^n\lambda_i\left[\langle v_0^\ast,
v_i-v_0\rangle+\phi(v_i)-\phi(v_0)+\intom j^0(x,\hat v_0(x);\hat
v_i(x)-\hat v_0(x))\; dx\right]<0.
\end{eqnarray*}
which is a contradiction. This implies that $F$ is a KKM mapping.

\item {\sc Step 2.} $F(v)\subseteq G(v)$ for all $v\in K\cap
D(\phi)$.

For a given $v\in K\cap D(\phi)$, let $u\in F(v)$. Then, there
exists $u^\ast \in A(u)$ such that
$$
\langle u^\ast, v-u\rangle+\phi(v)-\phi(u)+\intom j^0(x,\hat
u(x);\hat v(x)-\hat u(x))\; dx\geq 0.
$$
Since $A$ is monotone, we have
$$
\langle v^\ast-u^\ast, v-u\rangle\geq 0,\quad \forall v^\ast \in
A(v).
$$
It follows from the last two relations that
$$
\langle v^\ast, v-u\rangle+\phi(v)-\phi(u)+\intom j^0(x,\hat
u(x);\hat v(x)-\hat u(x))\; dx\geq 0,\quad\forall v^\ast A(v)
$$
which may be equivalently rewritten
$$
\sup\limits_{v^\ast\in A(v)} \langle v^\ast,u-v\rangle\leq
\phi(v)-\phi(v)+\intom j^0(x,\hat u(x);\hat v(x)-\hat u(x))\; dx
$$
and so $u\in G(v)$. This implies that $G$ is also a KKM mapping.

\item {\sc Step 3.} $G(v)$ is weakly closed for each $v\in K\cap
D(\phi)$.

Let $\{u_n\}\subset G(v)$ be a sequence which converges weakly to
$u$ as $n\ri \infty$. We must prove that $u\in G(v)$. Since
$u_n\in G(v)$ for all $n\geq 1$ and $\phi$ is weakly lower
semicontinuous, applying Lemma \ref{lema2}, for each $v^\ast\in
A(v)$ we have
\begin{eqnarray*}
0&\leq& \limsup\limits_{n\ri\infty}\left[ \langle v^\ast,
v-u_n\rangle+\phi(v)-\phi(u_n)+\intom j^0(x,\hat u_n(x);\hat
v(x)-\hat u_n(x))\; dx\right]\\
&\leq& \lim\limits_{n\ri\infty} \langle v^\ast,
v-u_n\rangle+\phi(v)-\liminf\limits_{n\ri\infty}\phi(u_n)
+\limsup\limits_{n\ri\infty}\intom
j^0(x,\hat u_n(x);\hat v(x)-\hat u_n(x))\;dx\\
&\leq& \langle v^\ast, v-u\rangle+\phi(v)-\phi(u)+\intom
j^0(x,\hat u(x);\hat v(x)-\hat u(x))\; dx.
\end{eqnarray*}
This is equivalent to $u\in G(v)$.

\item {\sc Step 4.} $G(v)$ is weakly compact for all $v\in K\cap
D(\phi)$.

Indeed, since $K$ is bounded, closed and convex, we know that $K$
is weakly compact, and so $G(v)$ is weakly compact for each $v\in
K\cap D(\phi)$, as it is a weakly closed subset of an weakly
compact set.
\end{description}
Therefore conditions of Lemma \ref{lema1} are satisfied in the
weak topology. It follows that
$$
\bigcap\limits_{v\in K\cap D (\phi)} G(v)\neq\emptyset.
$$
This yields an element $u\in K\cap D(\phi)$ such that, for any $v
\in K\cap D(\phi)$
$$
\sup\limits_{v^\ast\in A(v)}\langle v^\ast,
v-u\rangle\leq\phi(v)-\phi(u)+\intom j^0(x,\hat u(x);\hat
v(x)-\hat u(x))\; dx.
$$

This inequality is trivially satisfied for any $v\not\in D(\phi)$
which means that the inequality problem \eq{eq5} has at least one
solution. Theorem \ref{Th1} enables us to claim that inequality
problem \eq{eq4} also possesses a solution.
\endproof
\medskip

We shall derive below some consequences for problems over
unbounded sets. Without loss of generality we may assume that
$0\in K\cap D(\phi)$ and let us consider the sets $K_n=\{u\in K|\
\|u\|\leq n\}$ for $n=1,2,\ldots$

If $K$ is nonempty, unbounded, closed and convex subset of $X$ and
$A:K\ri 2^{X^\ast}$ is monotone and lower hemicontinuous, then by
Theorem \ref{Th2}, there exists $u_n\in K_n$ and $u_n^\ast\in
A(u_n)$ such that
\neweq{eq16}
\langle u_n^\ast, v-u_n\rangle+\phi(v)-\phi(u_n)+\intom j^0(x,\hat
u_n(x);\hat v(x)-\hat u_n(x))\; dx\geq 0,\quad \forall v\in K_n,
\endeq
for every $n\geq 1$.

\begin{teo}\label{Th3}
Assume that the same hypotheses as in Theorem \ref{Th2} hold
without the assumption of boundedness of $K$ and let $u_n\in K_n$
and $u_n^\ast\in A(u_n)$ be two sequences  such that \eq{eq16} is
satisfied for every $n\geq 1$. Then each of the following
condition is sufficient for the problem \eq{eq4} to possess a
solution:
\begin{description}
\item (H1) There exists a positive integer $n_0$ such that
$\|u_{n_0}\|<n_0$;

\item (H2) There exists a positive integer $n_0$ such that
$$
\langle u^\ast_{n_0},-u_{n_0}\rangle+\phi(0)-\phi(u_{n_0})+\intom
j^0(x,\hat u_{n_0}(x);-\hat u_{n_0}(x))\; dx\leq 0;
$$

\item (H3) There exists $u_0\in K\cap D(\phi)$ and $q\geq p$ such
that for any sequence $\{w_n\}\subset K$ with $\|w_n\|\ri \infty$
as $n\ri \infty$
$$
\frac{\langle w_n^\ast, w_n-u_0\rangle}{\|w_n\|^q}\ri \infty,\quad
\mbox{ as } n\ri\infty
$$
for every $w_n^\ast\in A(w_n)$.
\end{description}
\end{teo}

\proof Let $v\in K$ be arbitrary fixed.

\begin{itemize}
\item Assume (H1) holds and take $t>0$ small enough such that
$w=u_{n_0}+t(v-u_{n_0})$ satisfies $w\in K_{n_0}$ (it suffices to
take $t=1$ if $v=u_{n_0}$ and $t<(n_0-\|u_{n_0}\|)/\|v-u_{n_0}\|$
otherwise). By \eq{eq16} we have
\begin{eqnarray*}
0&\leq& \langle
u^\ast_{n_0},w-u_{n_0}\rangle+\phi(w)-\phi(u_{n_0})+\intom
j^0(x,\hat
u_{n_0}(x);\hat w(x)-\hat u_{n_0}(x))\; dx\\
&\leq& t\left[\langle
u^\ast_{n_0},v-u_{n_0}\rangle+\phi(v)-\phi(u_{n_0})+\intom
j^0(x,\hat u_{n_0}(x);\hat v(x)-\hat u_{n_0}(x))\; dx\right].
\end{eqnarray*}
Dividing by $t$ the last relation we observe that $u_{n_0}$ is a
solution of \eq{eq4}.

\item  Now, let us assume that (H2) is fulfilled. In this case,
some $t\in(0,1)$ can be found such that $tv\in K_{n_0}$. Taking
\eq{eq16} into account
\begin{eqnarray*}
0&\leq& \langle
u^\ast_{n_0},tv-u_{n_0}\rangle+\phi(tv)-\phi(u_{n_0})+\intom
j^0(x,\hat
u_{n_0}(x);t\hat v(x)-\hat u_{n_0}(x))\; dx\\
&=& \langle
u^\ast_{n_0},t(v-u_{n_0})+(1-t)(-u_{n_0})\rangle+
\phi(tv+(1-t)0)-\phi(u_{n_0})\\
&&+\intom j^0\left(x,\hat
u_{n_0}(x);t(\hat v(x)-\hat u_{n_0}(x))
+(1-t)(-\hat u_{n_0}(x)\right))\; dx\\
&\leq& t\left[\langle
u^\ast_{n_0},v-u_{n_0}\rangle+\phi(v)-\phi(u_{n_0})+\intom
j^0(x,\hat
u_{n_0}(x);\hat v(x)-\hat u_{n_0}(x))\; dx\right]\\
&&+(1-t)\left[ \langle
u^\ast_{n_0},-u_{n_0}\rangle+\phi(0)-\phi(u_{n_0})+\intom
j^0(x,\hat
u_{n_0}(x);-\hat u_{n_0}(x))\; dx\right]\\
&\leq& t\left[\langle
u^\ast_{n_0},v-u_{n_0}\rangle+\phi(v)-\phi(u_{n_0})+\intom
j^0(x,\hat u_{n_0}(x);\hat v(x)-\hat u_{n_0}(x))\; dx\right].
\end{eqnarray*}
Dividing again by $t$ the conclusion follows.

\item Assuming that (H3) holds we observe that there exists
$n_0>0$ such that $u_0\in K_n$ for all $n\geq n_0$. We claim that
the sequence $\{u_n\}$ is bounded. Suppose by contradiction that
up to a subsequence $\|u_n\|\ri \infty$. Since $w_n=u_n/\|u_n\|$
is bounded, passing eventually to a subsequence (still denoted
$w_n$ for the sake of simplicity), we may assume that $w_n\weak
w$. The function $\phi$ being convex and lower semicontinuous, it
is bounded from below by an affine and continuous function (see
e.g. \cite{bre} Prop. I.9), which means that for some $f\in
X^\ast$ and some $\alpha \in \RR$ we have
$$
\langle f,u\rangle+\alpha\leq\phi(u),\quad\forall u\in X.
$$
This leads to
\neweq{eq17}
-\phi(u)\leq \|f\|\cdot\|u\|-\alpha,\quad\forall u\in X.
\endeq
On the other hand, for any $y,h\in \RR^k$ there exists $z\in \de
j(x,y)$ (see \cite{Cla} Prop. 2.1.2) such that
$$
j^0(x,y;h)=z\cdot h=\max\{Z\cdot h|\ Z\in\de j(x,y)\}.
$$
It follows from \eq{eq1} that
$$
\left|j^0(x,\hat u(x);\hat v(x))\right|\leq C\left(1+|\hat
u(x)|^{p-1}\right)|\hat v(x)|
$$
and thus
\neweq{eq18}
\left|\intom j^0(x,\hat u(x);\hat v(x))\; dx\right|\leq
C\left(\|\hat v(x)\|_1+\|\hat u(x)\|_p^{p-1}\|\hat
v(x)\|_p\right)\leq C_1\|v\|+C_2\|u\|^{p-1}\|v\|
\endeq
for some suitable constants $C_1,C_2>0$. Relations \eq{eq16} ,
\eq{eq17}, \eq{eq18} show that
\begin{eqnarray*}
\langle u^\ast_n,u_n-u_0\rangle &\leq& \phi(u_0)-\phi(u_n)+\intom
j^0(x,\hat u_n(x);\hat u_0(x)-\hat u_n(x))\; dx\\
&\leq& \phi(u_0)+\|f\|\cdot \|u_n\|-\alpha
+C_1\|u_n-u_0\|+C_2\|u_n\|^{p-1}\|u_n-u_0\|.
\end{eqnarray*}
Thus
$$
\frac{\langle u_n^\ast, u_n-u_0\rangle}{\|u_n\|^q}\leq
\frac{\phi(u_0)-\alpha}{\|u_n\|^q}+\frac{\|f\|}{\|u_n\|^{q-1}}+
C_1\left\|\frac{w_n}{\|u_n\|^{q-1}}-\frac{u_0}{\|u_n\|^q}\right\|+
C_2\left\|\frac{w_n}{\|u_n\|^{q-p}}-\frac{u_0}{\|u_n\|^{q-p+1}}\right\|
$$
and passing to the limit as $n\ri\infty$ we reach a contradiction,
since $1<p\leq q$.

Since $\{u_n\}$ is bounded, a $n_0\geq 1$ can be found  such that
$\|u_{n_0}\|<n_0$ and by ($H_1$) $u_{n_0}$ solves \eq{eq4}.
\end{itemize}
\endproof

\end{document}